\documentclass[b5paper,twoside]{article}
\usepackage{amsmath, amsthm, amssymb}
\usepackage[dvips]{graphicx}
\def\titlePS#1{\noindent{{\LARGE\sffamily #1}} \vspace{3mm}}
\def\shorttitlePS#1{\markboth{Prague Stochastics 2006}{#1}}
\def\amsPS#1{\noindent{\it MSC 2000:} #1}
\def\keywPS#1{\noindent{\it Key words:} #1 \\}

\def\thanksPS#1{\footnotetext{#1}}
\def\authorsPS{}
\newcommand\authPS[8]{\ifnum\authPSc<1\def\authPSc{1}\else, \fi{\large\sffamily #1 #2}%
\edef\authorsPS{\authorsPS \par\vskip 1mm\noindent #1 #2:\hskip 5mm #3, #4, #5, #6, #7, #8}}
\newenvironment*{abstractPS}{\vspace{3mm}\noindent{\it Abstract:} }{\vspace{3mm}}
\newtheorem{theorem}{Theorem}[section]

\theoremstyle{definition}

\theoremstyle{remark}
\newtheorem{remark}[theorem]{Remark}

\def\authPSc{0}
\begin{document}

% !!!!!!!!!!!!!!!!!!!!!!!!!!!!!!!!!!!!!!!!!!!!!!!!!!!!!!!!!!!!!!!!!!!!!!!
% !!!                   START HERE                                   !!!
% !!!!!!!!!!!!!!!!!!!!!!!!!!!!!!!!!!!!!!!!!!!!!!!!!!!!!!!!!!!!!!!!!!!!!!!

%% !!  FULL TITLE OF PAPER
% example :
\titlePS{A Comparison of Information Concerning the Regression Parameter in The Accelerated
Failure Time Model under   Current Duration and Length Biased
Sampling: Does it Pay to be Patient?}

%% !! SHORT TITLE OF PAPER
%% To appear in headings, without mathematics please
\shorttitlePS{Current duration versus length biased sampling}

%% !!  ACKNOWLEDGMENT, SUPPORT ETC.
%%
\thanksPS{Acknowledgements.}

\noindent
%%   Author's data (copy for more authors)
%% !!  AUTHOR{forename}{family name}{institution}{street}{city}{ZIP code}{country}{e-mail}%
% example:

\authPS{Bert} {van Es} {University of Amsterdam, Korteweg-de Vries
Institute for Mathematics } {Plantage Muidergracht 24} {Amsterdam}
{1018 TV} {The Netherlands} {vanes @science.uva.nl}
\authPS{Chris
A.J. }{Klaassen}{University of Amsterdam, Korteweg-de Vries
Institute for Mathematics } {Plantage Muidergracht 24} {Amsterdam}
{1018 TV} {The Netherlands} {chrisk@science.uva.nl}
\authPS{Philip
J.}{Mokveld}{University of Amsterdam, Korteweg-de Vries Institute
for Mathematics } {Plantage Muidergracht 24} {Amsterdam} {1018 TV}
{The Netherlands}
%%%%%%%%%%%%%%%%%%%%%%%%%%%%%%%%%%
%  personal definitions of symbols like for example:
%  \def\Cb{{\mathcal C}_b}

%  are possible, if you need any, please, put them here
%%%%%%%%%%%%%%%%%%%%%%%%%%%%%%%%%%

\begin{abstractPS}
Longitudinal observations are sometimes costly  or not available.
Cross sectional sampling can be an alternative. Observations are
drawn  then at a specific point in time from a population of
durations whose distributions satisfy a {\em core model}.
Subsequently, one has a choice. One may process the data
immediately, obtaining so called current duration data. Or one
waits until the sampled durations are known completely   obtaining
the full  durations via length biased sampling. We compare the
Fisher information for the Euclidean parameter corresponding to an
Accelerated Failure Time core model when the observations are
obtained by either
 current duration  or length biased sampling.
\end{abstractPS}

%  AMS classification see http://www.ams.org/msc/index.html#search
\amsPS{62N02, 62D05}

%  Key words
\keywPS{Survival analysis, semiparametric statistics, cross
sectional sampling}

%% END OF TITLE INFORMATION

%% START OF MAIN BODY

\section{Current duration and length biased sampling from the AFT model}

Two often used models in survival analysis based on longitudinal
data are the Cox Proportional Hazards model (PH) and the
Accelerated Failure Time (AFT) model. These two semiparametric
models both have appealing interpretations and their properties
are well understood. For instance information bounds and efficient
estimators of the Euclidean regression parameter are available for
both models.

In  situations where longitudinal observations are costly, or not
available, one has to resort to technically more complicated but
less costly sampling schemes, like cross sectional sampling. In a
medical setting this means that instead of following a certain
number of patients in time one selects the durations of the
disease of a group of patients sampled at a specific point in
time, obtaining a so called {\em cross sectional sample}. One then
has a choice. Either one uses the data at hand at the time of
sampling, i.e. the durations up to the present, obtaining so
called {\em current duration data}, or one decides to wait until
the full durations for the sampled patients are known. Because
longer durations turn out to be sampled more frequently than
shorter ones, the second type of sampling is known as {\em length
biased sampling}.

Let us compare the  two cross sectional sampling regimes. Current
duration sampling will only require knowledge of the duration up
to the present and is thus very cheap in this sense. Length biased
sampling requires the time needed to observe the full durations of
the diseases of the patients that have been sampled and is thus
more costly than current duration sampling.

We will assume that we sample from a population of durations that
satisfy a semiparametric {\em core model}.  By comparing information
bounds for the Euclidean parameter under the two cross sectional
sampling schemes we will investigate the gain in efficiency in being
patient.

Our comparison below is based on results for current duration and
length biased sampling for the AFT core model in these situations,
presented in Mokveld (2006). Similar results for the PH model do not
exist at present.  See also Van Es, Klaassen and Oudshoorn
(2000) for some general features of current duration sampling.

\subsection{The core AFT model}

We first introduce the AFT core model. Let $T$ denote a duration,
for instance the duration of the disease of an individual from a
homogeneous group of patients with a particular disease, and let
$W$ denote  a vector of covariates of dimension $k$ with density
$h$ with respect to a measure $\nu$. We do not assume knowledge of
$h$. Let $\theta\in \Theta$ denote an unknown $k$-vector of
regression parameters.

The semiparametric AFT model for the random vector $(T,W)$ is
given by
\begin{equation}\label{eq:1}
        T=e^{-\theta^T W}V,
\end{equation}
where $V$ is a nondegenerate random variable on $[0,\infty)$ with
unknown absolutely continuous distribution function $G_0$, with
density $g_0$ and hazard function $\lambda_0$, and where $V$ and
$W$ are independent. We consider estimation of $\theta$, treating
$g_0$ as a nuisance parameter .

From the model equation (\ref{eq:1}) we can derive the conditional
survival function $\bar{G}_\theta(t|w)$, the conditional density
$g_\theta(t|w)$ and the conditional hazard function
$\lambda_\theta(t|w)$ of $T$ given $W=w$. We get, for $t>0$,
\begin{align*}
&\bar{G}_\theta(t|w) =1-
G_\theta(t|w)=\bar{G}_0(e^{\theta^Tw}t),\\
&g_\theta(t|w) =e^{\theta^T w} g_0(e^{\theta^T w}t),\\
&\lambda_\theta(t|w) =e^{\theta^T w}\lambda_0(e^{\theta^T w}t).
\end{align*}
Note that given the value of the covariate vector the model is a
scale model. The function $\lambda_0$ serves as baseline hazard in
this scale model. Depending on the value of the scale $e^{\theta^T
w}$ on average the duration is decreased or increased.

Also note that taking logarithms in the model equation (\ref{eq:1})
we get
$$
        \ln T=-\theta^T W + \ln V,
$$
showing that the AFT model is actually a regression model for the
logarithm of the duration. However, differences are caused by
different natural assumptions on the distributions of $V$ in the AFT
model and the error $\ln V$ in the regression model.

\subsection{Current duration and length biased sampling}

Let us assume that we observe the durations and their covariates at
a {\em specific point in time, the present}. Let $D$ denote the
total length of a sampled duration and let $X$ denote the time from
onset until the present of a sampled duration.

\bigskip

For simplicity we first describe the sampling distributions in the
situation without covariates. If $f$ and $F$ are the density and
distribution function of the durations $T$ in the core model then
under suitable assumptions the densities of $D$ and $X$ equal
\begin{align}
&\label{lb} f_D(y) =\frac{y f(y)}{\mu} ,\\
&\label{cs} f_X(x)=\frac{{\bar F}(x)}{\mu},
\end{align}
where ${\bar F}(x)=1-F(x)$ and $\mu=\int_0^\infty uf(u)du$.  It
turns out that $X$ is in distribution equal to $DU$  with $U$
uniformly distributed on the unit interval and with $D$ and $U$
independent. Hence, while formula (\ref{lb}) follows from the length
bias in the sampling, formula (\ref{cs}) follows from the same length
bias in  selecting   the duration and  from multiplicative
censoring, since at the present we
only observe a fraction of the total duration!

The formulas (\ref{lb}) and (\ref{cs}) require suitable models for
the times of onset of the disease. In Van Es, Klaassen and
Oudshoorn (2000) and Mokveld (2006) two models for the times of
onset are described that give rise to the densities above.

One can follow a {\em direct approach}   where the random variable
$L$ denotes the time of onset and is uniformly distributed on the
interval $[-\tau,0]$. Subsequently one lets $\tau$ go to infinity.
The duration $T$ is assumed to be independent from $L$ and current
duration sampling takes place at time zero. A duration is sampled if
and only if $T\geq-L$ (random left truncation). The disease will
have lasted $X=-L$ at time zero and will last $D=T$ if we wait until
recovery. The distributions of $X$ and $D$ can be computed by
conditioning on $T\geq -L$.

Following Keiding (1991) one can also follow a {\em point process
approach} where patients get ill at the time points of a
stationary Poisson   process with constant intensity $\lambda$.
The durations of their disease are modelled as i.i.d random
variables $T$ that are independent from the Poisson  process and
cross sectional  sampling takes place at some fixed point in time.
By point process techniques one can show that $N$, the number of
durations that are sampled, has a Poisson distribution, and,
conditionally on $N=n$, the sampled times $X$ from onset  and full
durations $D$ are i.i.d. with the densities (\ref{lb}) and
(\ref{cs}).

\bigskip

In the regression setting with covariates we observe $n$ i.i.d.
realizations of  $(D,Z)$ or $(X,Z)$ of durations (in total or from
onset to present) and the sampled covariates. As mentioned above
we consider the case where the density $h$ of the covariate $W$ in
the core model is unknown.  Under the AFT model assumptions for
the core model, it turns out that given the covariate $Z$ the
distributions of both $D$ and $X$ belong to scale parameter
families, just as the distribution of the original durations  $T$
in the core model. In fact, they again follow an AFT model. The
difference with the core model is that now   the distribution of
$Z$, the observed covariate, depends on the Euclidean parameter
$\theta$. It does not depend on $g_0$!

For $x>0,\ y>0$, and $z\in \mathbb{R}^k$ we have for the total
duration $D$
\begin{align}
        &f_{D,Z}(y,z)=\frac{e^{\theta^T z}y g_0(e^{\theta^T z}y)h(z)}{E_{g_0} V E_h e^{-\theta^TW}}, \nonumber\\
        &\label{cov:1}f_{Z}(z)=\frac{e^{-\theta^T z} h(z)}{E_h e^{-\theta^T W}}, \\
        &f_{D|Z}(y|z)=\frac{ e^{2\theta^Tz}y
        g_0(e^{\theta^Tz}y)}{E_{g_0}V},\nonumber
\end{align}
and for the duration from onset to present $X$
\begin{align}
        &f_{X,Z}(x,z)=\frac{\bar{G_0}(e^{\theta^T z}x)h(z)}{E_{g_0} V E_h e^{-\theta^TW}} , \nonumber \\
        &\label{cov:2}f_{Z}(z)=\frac{e^{-\theta^T z} h(z)}{E_h e^{-\theta^T W}},  \\
        &f_{X|Z}(x|z)=\frac{e^{\theta^T z} \bar{G}_0(e^{\theta^T z}x)}{E_{g_0} V}. \nonumber
\end{align}
These formulas hold under the direct approach or the point process
approach for the times of onset described above. See Van Es,
Klaassen and Oudshoorn (2000) or Mokveld (2006) for details.

\section{A comparison of information bounds}

We will present information bounds for estimation of the Euclidean
parameter $\theta$ for cross sectional sampling from a core AFT
model as derived in Mokveld (2006).  Throughout, when we mention
information we mean information contained in one observation.

As above, primarily we consider the case where the covariate
distribution is unknown. See Remark \ref{hknown} for the case where
this distribution is known.

\subsection{Current duration and length biased sampling}

The covariance matrix of the sampled covariates appears in all
information matrices below. It equals
$$
\Sigma_Z=E (Z-E Z)(Z-E Z)^T.
$$
Note that this matrix  depends on $\theta$ through the distribution
of $Z$.

Let us first define the {\em Fisher information for scale} $I_s(f)$
for a   density $f$
\begin{equation}\label{scale}
I_s(f)=\int \Big(1+x\,\frac{f'(x)}{f(x)}\Big)^2f(x)dx.
\end{equation}
With $\mu=\int xg_0(x)dx= \int {\bar G_0}(x)dx$,  $f_1(x)$ equal
to $xg_0(x)/\mu$ and $f_2(x)$ equal to ${\bar G_0(x)}/\mu$, it is
shown in Mokveld (2006) that efficient estimators of $\theta$ can
be constructed and that the information bounds are equal to
$$
\Sigma_Z\, I_s(f_1)
$$
in the situation of length biased sampling where the full durations
are observed, and to
$$
\Sigma_Z\, I_s(f_2)
$$
in the situation of current duration sampling where the durations
from onset to present are observed. Rewriting  $I_s(f_1)$ and
$I_s(f_2)$ in terms of $g_0$ we get
$$
I_s(f_1)=\int
\Big(2+x\,\frac{g_0'(x)}{g_0(x)}\Big)^2\frac{xg_0(x)}{\mu}dx
$$
and
$$
I_s(f_2)=\int \Big(1-x\,\frac{g_0'(x)}{{\bar
G_0}(x)}\Big)^2\frac{{\bar G}_0(x)}{\mu}dx.
$$

\begin{remark}\label{hknown}
Let us consider the model where the covariate distribution in the
core model is known. Then (\ref{cov:1}) and (\ref{cov:2}) show
that the distribution of the covariates $Z$ in the sample is the
same for current duration and length biased sampling,   that it
does not depend on $g_0$, and that the Fisher information matrix
in one observation for $\theta$ based on the covariates in the
sample alone is   equal to $\Sigma_Z$. Under suitable assumptions
$\theta$ can be estimated $\sqrt{n}$-consistently from the
covariates alone by for instance the maximum likelihood estimator.

 The information for $\theta$ based on durations and
covariates now equals
$$
\Sigma_Z(I_s(f_1)+1)
$$
in the situation of length biased sampling where the full durations
are observed, and to
$$
\Sigma_Z(I_s(f_2)+1)
$$
in the situation of current duration sampling where the durations
from onset to present are observed. These are obviously larger than
in the situation where the covariate distribution is unknown.

Note also that, using both durations and covariates in the sample,
the semiparametric information for $\theta$, with $g_0$ as
nuisance parameter, under the two sampling schemes,  is larger
than $\Sigma_Z$, the information based on the covariates alone.
\end{remark}

\subsection{A comparison}

The results in this section show that it pays to be patient.

\begin{theorem}\label{thm:1} Let $g$ be an absolutely continuous density on
$(0,\infty)$ with derivative $g'$ a.e. and let $\mu=\int xg(x)dx<
\infty$. Let $f_1(x)$ be equal to $xg(x)/\mu$ and let $f_2(x)$ be
equal to ${\bar G}(x)/\mu$. If $I_s(f_2)$ and $I_s(f_1)$ are
finite then
\begin{equation}\label{inequality}
 I_s(f_2)<I_s(f_1)
\end{equation}
holds.
\end{theorem}

\begin{proof} Note that $f_1$ is the density of $Y_1=e^{\theta^TZ}X$
and that $f_2$ is the density of $Y_2=e^{\theta^TZ}D$. Since
$X=UD$, with $U$ independent of $D$ and uniformly distributed on
the unit interval, we have
$$
P(Y_2\leq x) = \int_0^1 P\Big(Y_1\leq\frac{x}{u}\Big)du.
$$
So the relation between   $f_1$ and  $f_2$ can be expressed as
\begin{equation}\label{f1dens}
f_2(x)=\int_0^1 f_1\Big(\frac{x}{u}\Big)\frac{1}{u}du.
\end{equation}
By expanding the square  in (\ref{scale})  we see that the
inequality (\ref{inequality}) holds if and only if
\begin{equation}\label{ff}
\int x^2\Big(\frac{f_2' }{f_2}\Big)^2 (x)f_2(x)  dx <   \int
x^2\Big(\frac{f_1'}{f_1}\Big)^2 (x)f_1(x)  dx.
\end{equation}
Let $f_1$ vanish at $x_0$ and be differentiable at $x_0$ with
derivative $f_1'(x_0)$. Since $f_1$ is nonnegative Lebesgue a.e.,
we get $f_1'(x_0)=0$. Because an absolutely continuous function is
Lebesgue a.e. differentiable, this shows that $\{ x:
f_1(x)=0,f_1'(x)\not= 0\}$ is a Lebesgue null set. Consequently by
the Cauchy-Schwarz inequality we have
\begin{align*}
\Big(\int_0^1&f_1'\Big( \frac{x}{u}\Big)\frac{1}{u^2}du\Big)^2 =
\left(\int_0^1 \left\{ \frac{ f_1'\Big(
\frac{x}{u}\Big)\frac{1}{u^2}}{\sqrt{f_1\Big(
\frac{x}{u}\Big)\frac{1}{u}}} \right\}\left\{\sqrt{f_1\Big(
\frac{x}{u}\Big)\frac{1}{u}}\right\}du\right)^2\\
&\leq\int_0^1 \Big(\frac{f_1'}{f_1}\Big)^2 (\frac{x}{u})  f_1\Big(
\frac{x}{u}\Big)\frac{1}{u^3}du\int_0^1\frac{1}{u}f_1\Big(\frac{x}{u}\Big)
du.
\end{align*}
Hence by (\ref{f1dens}) we have
\begin{align*}
\int x^2&\Big(\frac{f_2'}{f_2}\Big)^2  (x) f_2(x) dx = \int x^2
\left\{\frac{ \int_0^1 f_1'\Big(\frac{x}{u}\Big)\frac{1}{u^2}du }
{\int_0^1 f_1\Big(\frac{x}{u}\Big)\frac{1}{u}du}\right\}^2
\int_0^1 f_1\Big(\frac{x}{u}\Big)\frac{1}{u}dudx\\
& \leq\int\int_0^1\frac{x^2}{u^3}\Big(\frac{ f_1' }{ f_1
}\Big)^2\Big( \frac{x}{u}\Big)f_1\Big( \frac{x}{u}\Big)dudy=
\int_0^1\int x^2\Big(\frac{ f_1' }{ f_1  }\Big)^2 (x)f_1(x)dxdu\\
& =\int x^2\Big(\frac{f_1'}{f_1}\Big)^2  (x) f_1(x) dx,
\end{align*}
which completes the proof of the inequality provided that we show that equality can not occur.

The fact that the inequality (\ref{inequality}) is strict can be seen as follows.
The Cauchy-Schwarz inequality holds with equality if and only if
$$
\frac{ f_1'\Big( \frac{x}{u}\Big)\frac{1}{u^2}}{\sqrt{f_1\Big(
\frac{x}{u}\Big)\frac{1}{u}}} =c\sqrt{f_1\Big(
\frac{x}{u}\Big)\frac{1}{u}},
$$
for some constant $c$ and for all $u\in[0,1]$. But for equality to
hold in (\ref{ff}) this last equality has to hold for all $x$. Now
writing $z=x/u$ this condition equals
$$
zf_1'(z)=cxf_1(z)
$$
for all $x>0$ and all $z>x$, which can obviously never hold.
\end{proof}

\bigskip

 Actually, this theorem is a consequence of a more general inequality for
Fisher information for scale for a product of random variables.
\begin{theorem}
Let $f$ be a density on $(0,\infty)$ that is absolutely continuous
with respect to Lebesgue measure with derivative $f'$, such that
$I_s(f)$, as defined by (\ref{scale}), is finite. If $G$ is an
arbitrary distribution function on $(0,\infty)$ and the density $h$
is defined by
$$h(x)=\int_{0}^\infty {1 \over u} f\Big({x \over u}\Big) dG(u)$$
then
$$ I_s(h) \leq I_s(f)$$
with equality iff $G$ is degenerate.
\end{theorem}

\begin{proof} Let $X$ be a random variable with density $f\,.$ The
random variable $\log X$ has density $\tilde f$ then with ${\tilde
f}(x) = e^x f(e^x)\,.$ One may verify that the Fisher information
$I_s(f)$ for scale of $f$ equals the Fisher information
$I_\ell({\tilde f})$ for location of $\tilde f \,.$ Furthermore,
$h$ is the density of the product of $X$ and a random variable
with distribution $G\,.$ Consequently, with $\tilde h$ defined by
${\tilde h}(z)=e^zh(e^z)$ and $\tilde G (z)$ defined by $G(e^z)$,
it suffices to prove that
$I_\ell({\tilde h}) \leq I_\ell({\tilde f})$ holds with equality
iff $\tilde G$ is degenerate. However,
this inequality follows by
Cauchy-Schwarz via
\begin{align*}
\int \Big(\frac{\tilde h'}{\tilde h}\Big)^2\tilde h
&= \int \frac{ \Big\{ \int \frac{\tilde f'}{\tilde f} (x-y) \sqrt{\tilde f(x-y)}
\sqrt{\tilde f(x-y)} d{\tilde G}(y)    \Big\}^2 }{\tilde h(x)}\, dx\\
&\leq
\int\int  \Big(\frac{\tilde f'}{\tilde f}\Big)^2(x-y) {\tilde f} (x-y)   d{\tilde G}(y)dx
=I_l(\tilde f),
\end{align*}
as has been noticed by H\'ajek and \v{S}id\'ak (1967) in their
Theorem I.2.3 on page 17.
\end{proof}

\subsubsection{Examples}

To get a feeling for the difference in information in the current
duration and length biased observations we consider two families
of densities for the nuisance parameter $g_0$, the Weibull
densities and the log logistic densities.

First we consider the Weibull densities. Let $g_0$ be a Weibull density
with parameter $\gamma >0$, i.e.
$$
g_0(t)=\gamma t^{\gamma -1}e^{-t^{\gamma}}, \quad t\geq 0.
$$
For these densities we have
\begin{align*}
&I_s(f_1)=\gamma(\gamma +1) ,\\
&I_s(f_2)=\gamma .\\
\end{align*}

\begin{figure}
\includegraphics[width=6cm]{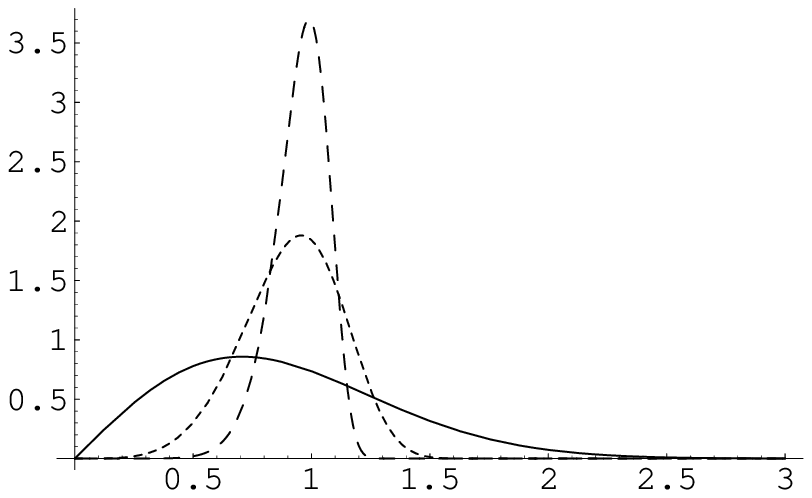}
\includegraphics[width=6cm]{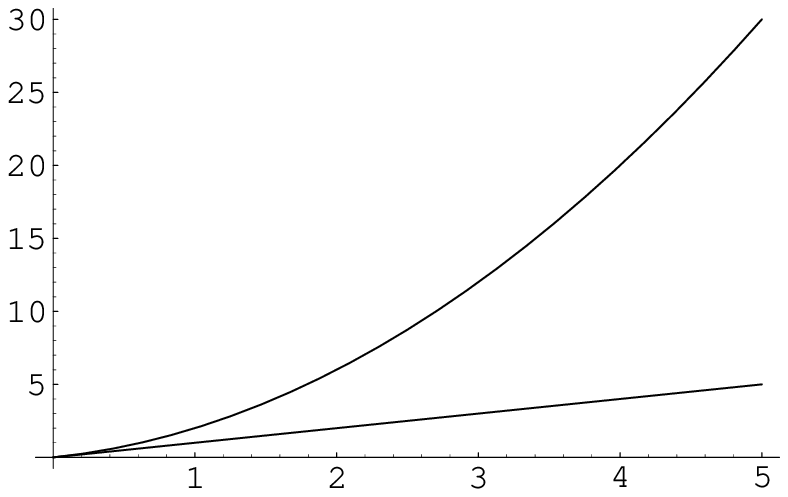}
\caption{Left: Weibull densities for $\gamma$ equal to 2 (solid), 5
(...) and 10 (- - -). Right: information under length biased and
current duration sampling (Weibull $g$) as a function of $\gamma$.}
\end{figure}

Next we consider log logistic densities $g_0$. Let $g_0$ be a log
logistic density with parameter $\gamma >1$, i.e.
$$
g_0(t)=\frac{\gamma t^{\gamma -1}}{(1+t^\gamma)^2}, \quad t\geq 0.
$$
For these densities we have
\begin{align*}
&I_s(f_1)=\frac{1}{3}(\gamma^2-1) ,\\
&I_s(f_2)=\frac{1}{2}(\gamma-1)  .\\
\end{align*}
\begin{figure}
\includegraphics[width=6cm]{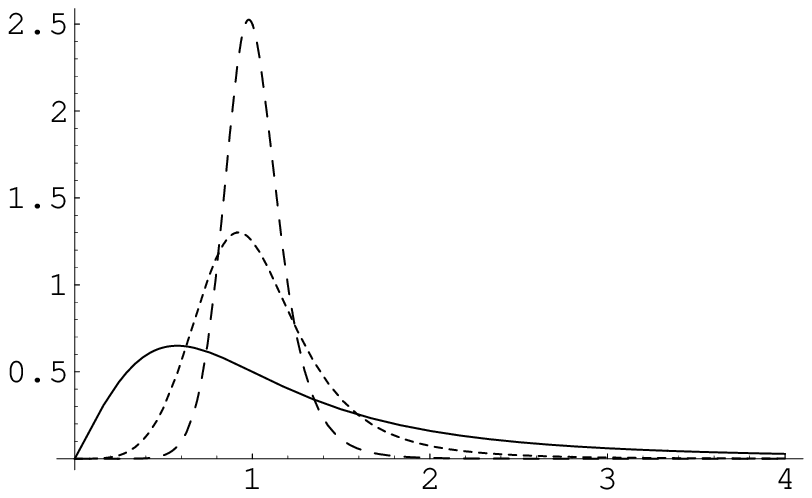}
\includegraphics[width=6cm]{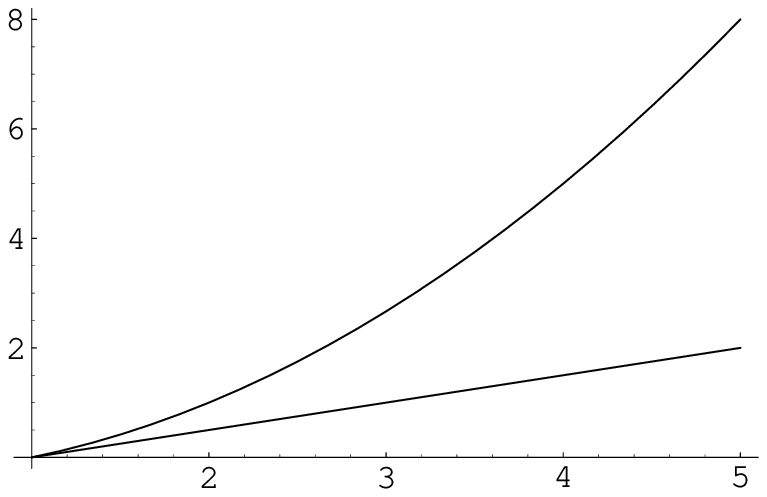}
\caption{Left: Log logistic densities for $\gamma$ equal to 2
(solid), 5 (...) and 10 (- - -). Right: information under length
biased and current duration sampling (log logistic $g$) as a
function of $\gamma$.}
\end{figure}

These two examples show  that the more concentrated the density
$g_0$ of the random variable $V$ in the model (\ref{eq:1}),
corresponding with high parameter values $\gamma$,
the higher the gain in being patient.

%% the bibliography is formatted using \bibligraphystyle{plain}

% !!!!!!!!!!!!!!!!!!!!!!!!!!!!!!!!!!!!!!!!!!!!!!!!!!!!!!!
% !!!                   END HERE                      !!!
% !!!!!!!!!!!!!!!!!!!!!!!!!!!!!!!!!!!!!!!!!!!!!!!!!!!!!!!
\authorsPS

\begin{thebibliography}{10}
% example:

\bibitem{eko:2000} B. van Es, C.~A.~J. Klaassen, and K. Oudshoorn.
\newblock Survival analysis under cross sectional sampling.
Prague Workshop on Perspectives in Modern Statistical Inference: Parametrics, Semi-parametrics, Non-parametrics
(1998).
\newblock {\em J. Statist. Plann. Inf.}, 91:295--312, 2000.

\bibitem{hajeksidak:1967}
J. H\'ajek and Z. \v{S}id\'ak (1967).
\newblock {\em Theory of Rank Tests.}
\newblock Academic Press, New
York-London; Academia Publishing House of the Czechoslovak Academy
of Sciences, Prague, 1967.


\bibitem{keid:1991} N. Keiding.
\newblock Age specific incidence and prevalence: a statistical perspective.
\newblock {\em J. Roy. Statist. Soc. Ser. A}, 154:371--412, 1991.

\bibitem{mok:2006}
Ph.~J. Mokveld.
\newblock {\em The Accelerated Failure Time Model under Cross
Sectional Sampling Schemes}.
\newblock Ph.D. Thesis, University of Amsterdam, in preparation.



\end{thebibliography}
\end{document}